# COMPOSITE QUANTILE REGRESSION AND THE ORACLE MODEL SELECTION THEORY[1]

By Hui Zou and Ming Yuan

*University of Minnesota and Georgia Institute of Technology*

Coefficient estimation and variable selection in multiple linear regression is routinely done in the (penalized) least squares (LS) framework. The concept of model selection oracle introduced by Fan and Li [*J. Amer. Statist. Assoc.* **96** (2001) 1348–1360] characterizes the optimal behavior of a model selection procedure. However, the least-squares oracle theory breaks down if the error variance is infinite. In the current paper we propose a new regression method called *composite quantile regression* (CQR). We show that the oracle model selection theory using the CQR oracle works beautifully even when the error variance is infinite. We develop a new oracular procedure to achieve the optimal properties of the CQR oracle. When the error variance is finite, CQR still enjoys great advantages in terms of estimation efficiency. We show that the relative efficiency of CQR compared to the least squares is greater than 70% regardless the error distribution. Moreover, CQR could be much more efficient and sometimes arbitrarily more efficient than the least squares. The same conclusions hold when comparing a CQR-oracular estimator with a LS-oracular estimator.

## 1. Introduction and motivation.

1.1. *Background.* In recent years, various techniques have been developed for simultaneous variable selection and coefficient estimation in multiple linear regression. Notable methods include the nonnegative garrote [Breiman (1995)], the lasso [Tibshirani (1996)] and the SCAD [Fan and Li (2001)]. Fan and Li (2006) gave a comprehensive overview of recent advances in variable selection.

Received March 2007; revised May 2007.
[1]Supported in part by NSF Grants DMS-07-06733 and DMS-07-60724.

Fan and Li (2001) introduced the concept of model selection oracle to guide the construction of optimal model selection procedures. To elaborate, consider the following linear model

$$y = \sum_{j=1}^{p} x_j \beta_j^* + \varepsilon. \tag{1.1}$$

Without loss of generality we center the predictors. Denote $\mathcal{A} = \{j : \beta_j^* \neq 0\}$. The problem of variable selection and coefficient estimation is to identify the unknown set $\mathcal{A}$ and estimate the corresponding coefficients $\boldsymbol{\beta}_{\mathcal{A}}^*$, using $n$ independent samples generated from the model (1.1). To understand the optimality of variable selection and coefficient estimation, Fan and Li (2001) suggested considering the oracle who knows the true subset $\mathcal{A}$. The oracle would only need to estimate $\boldsymbol{\beta}_{\mathcal{A}}^*$ and set $\boldsymbol{\beta}_{\mathcal{A}^c}^* = 0$. It is worth emphasizing here that although the oracle knows the true subset model, the error distribution remains unknown. Fan and Li (2001) considered the oracle estimator which estimates $\boldsymbol{\beta}_{\mathcal{A}}^*$ by least squares, which we shall refer to as the *LS-oracle*. Denote by $\mathbf{X}$ the design matrix and assume that $\lim_{n \to \infty} \frac{1}{n} \mathbf{X}^T \mathbf{X} = \mathbf{C}$, where $\mathbf{C}$ is a $p \times p$ positive definite matrix. Write $\mathbf{C}_{\mathcal{A}\mathcal{A}}$ the sub-matrix of $\mathbf{C}$ with both row and column indices in $\mathcal{A}$. We have

$$\sqrt{n}(\widehat{\boldsymbol{\beta}}^{\mathrm{LS}}(oracle)_{\mathcal{A}} - \boldsymbol{\beta}_{\mathcal{A}}^*) \to_d N(0, \sigma^2 \mathbf{C}_{\mathcal{A}\mathcal{A}}^{-1}), \tag{1.2}$$

where $\sigma^2$ is the variance of $\varepsilon$.

Note that the oracle "estimator" is not a legitimate estimator because it uses the information of $\mathcal{A}$, which is unavailable in practice. Nevertheless, the oracular model selection theory provides a golden standard for evaluating variable selection and coefficient estimation procedures. Following Fan and Li (2001), we say a variable selection and coefficient estimation procedure $\eta$ is a *LS-oracular estimator*, if $\widehat{\boldsymbol{\beta}}(\eta)$ (asymptotically) has the following oracle properties:

- Consistent selection: $\Pr(\{j : \hat{\boldsymbol{\beta}}(\eta)_j \neq 0\} = \mathcal{A}) \to 1$.
- Efficient estimation: $\sqrt{n}(\widehat{\boldsymbol{\beta}}(\eta)_{\mathcal{A}} - \boldsymbol{\beta}_{\mathcal{A}}^*) \to_d N(0, \sigma^2 \mathbf{C}_{\mathcal{A}\mathcal{A}}^{-1})$.

Thus $\eta$ works as well as the LS-oracle. Fan and Li (2001) showed that the SCAD indeed attains the oracle properties. Zou (2006) later demonstrated that the adaptive lasso also enjoys the oracle properties.

1.2. *Issues with the oracle.* Ideally, if one knows the error distribution, then the best oracle is the maximum likelihood estimate knowing the true underlying sparse model. However, in the linear regression problems, the error distribution (hence the likelihood model) is typically unknown. Hence we can only consider a practical oracle procedure. Although Fan and Li (2001) treated the (penalized) least squares as a special case in a general (penalized) likelihood based framework when the noise follows a normal distribution, we should note that the LS-oracle model selection theory does not need the normal error assumption, as long as the error distribution has a finite variance. However, the finite variance assumption is crucial for the oracle model selection theory based on the least squares. The reason is simple. If the error variance is infinite, $\widehat{\boldsymbol{\beta}}^{\text{LS}}(oracle)$ is no longer a root-$n$ consistent estimator and can not serve as an ideal method for variable selection and coefficient estimation. Consequently, estimators that are shown to enjoy the oracle properties, such as the SCAD or the adaptive lasso, will also lose their root-$n$ consistency. Such situations occur, for example, when the error follows a Cauchy distribution. On the other hand, model selection is about discovering the sparse structure in the relation between the response and predictors; thus model selection is still a legitimate and interesting problem even when the error variance is infinite. See Example 5 in Section 5.

The limitation of the LS-oracle motives our work here. We wish to find an alternative oracle that can overcome the breakdown issue of the LS oracle. There are several important considerations when designing a new oracle estimator. Let $\widehat{\boldsymbol{\beta}}^{\text{new}}(oracle)$ be a new oracle estimator. Firstly, the new oracle estimator should be root-$n$ consistent and enjoy asymptotic normality even when the LS-oracle fails to do so. Secondly, we are interested in the relative efficiency of the new oracle estimator $\widehat{\boldsymbol{\beta}}^{\text{new}}(oracle)$ with respect to $\widehat{\boldsymbol{\beta}}^{\text{LS}}(oracle)$ when $\sigma^2 < \infty$. Since $\widehat{\boldsymbol{\beta}}^{\text{LS}}(oracle)$ is of full efficiency when the error follows a normal distribution, it is impossible to have an oracle that is universally more efficient than the LS-oracle. However, it would be very nice to have the relative efficiency of $\widehat{\boldsymbol{\beta}}^{\text{new}}(oracle)$ with respect to $\widehat{\boldsymbol{\beta}}^{\text{LS}}(oracle)$ be bounded from below. This will prevent severe loss of statistical efficiency even in the worst scenario. Furthermore, we would like to see that $\widehat{\boldsymbol{\beta}}^{\text{new}}(oracle)$ can be significantly more efficient than $\widehat{\boldsymbol{\beta}}^{\text{LS}}(oracle)$ for commonly used nonnormal error distributions. Finally, the oracle estimator needs to be attainable in the sense that we have an estimating procedure that can mimic $\widehat{\boldsymbol{\beta}}^{\text{new}}(oracle)$, like the SCAD mimics $\widehat{\boldsymbol{\beta}}^{\text{new}}(LS)$.

It is not a trivial task to find an oracle estimator that satisfies all the above properties. For instance, the least absolute value regression is an obvious alternative to the least squares. Even for Cauchy-distributed errors, the least absolute value regression estimator still enjoys the asymptotic normality. The oracle estimator by the least absolute value regression is also attainable by the SCAD [see Fan and Li (2001), page 1357]. However, the relative efficiency of the least absolute value regression can be arbitrarily small when compared with the least squares. Therefore, we do not consider it as a safe alternative to the least squares.

1.3. *Our contributions.* In this work we introduce a new regression method called *composite quantile regression* (CQR) that can be used to construct an oracle estimator possessing all of the aforementioned properties. We define the CQR in Section 2. In Section 3 we study the asymptotic relative efficiency of the CQR-oracle with respect to the LS-oracle. A universal lower bound is derived, which shows that the relative efficiency is always larger than 70%. Moreover, we show by several concrete examples that the CQR-oracle can be much more efficient than the LS-oracle. In Section 4 we propose to use the adaptive lasso penalty to construct the adaptively penalized CQR estimator which is shown to achieve the performance of the CQR-oracle if the penalization parameter is appropriately chosen. Simulation results are presented in Section 5. The technical proofs are presented in Section 6. Section 7 contains a few concluding remarks.

**2. Composite quantile regression.** To motivate CQR, let us briefly review the quantile regression method [Koenker (2005)]. Note that the conditional $100\tau\%$ quantile of $y|\mathbf{x}$ is

$$\sum_{j=1}^{p} x_{ij}\beta_j^* + b_\tau^*$$

where $b_\tau^*$ is the $100\tau\%$ quantile of $\varepsilon$. For brevity, we shall assume that the density function of $\varepsilon$ is nonvanishing everywhere. Therefore $b_\tau^*$ is uniquely defined for any $0 < \tau < 1$. Quantile regression estimates $\boldsymbol{\beta}^*$ by solving

$$(2.1) \qquad (\hat{b}_\tau, \widehat{\boldsymbol{\beta}}^{\mathrm{QR}_\tau}) = \arg\min_{b,\boldsymbol{\beta}} \sum_{i=1}^{n} \rho_\tau\left(y_i - b - \sum_{j=1}^{p} x_{ij}\beta_j\right),$$

where $\rho_\tau(t) = \tau t_+ + (1-\tau)t_-$ is the so-called check function where subscripts $+$ and $-$ stand for the positive and negative parts, respectively. Quantile regression has been widely used in various areas such as economics [Koenker and Hallock (2001)] and survival analysis [Koenker and Geling (2001)] among others. It is well known that under mild regularity conditions [Koenker (2005)],

$$(2.2) \qquad \sqrt{n}(\widehat{\boldsymbol{\beta}}^{\mathrm{QR}_\tau} - \boldsymbol{\beta}^*) \to_d N\left(0, \frac{\tau(1-\tau)}{f^2(b_\tau^*)}\mathbf{C}^{-1}\right).$$

Quantile regression can be more efficient than the least squares estimator. In particular, if $\varepsilon$ follows a double-exponential distribution, $\widehat{\boldsymbol{\beta}}^{\mathrm{QR}_{0.5}}$ is the most efficient estimator. $\widehat{\boldsymbol{\beta}}^{\mathrm{QR}_\tau}$ does not require $\sigma^2 < \infty$ in order to enjoy the root-$n$ consistency and asymptotic normality, as opposed to the LS estimator. However, the relative efficiency of the quantile regression estimator with respect to the LS estimator can be arbitrarily small.

To further improve upon the usual quantile regression, we propose to simultaneously consider multiple quantile regression models. Note that the regression coefficients are the same across different quantile regression models. We shall demonstrate that by combining the strength across multiple quantile regression models, we can derive a good estimator that satisfies the desired properties discussed in the introduction.

Denote $0 < \tau_1 < \tau_2 < \cdots < \tau_K < 1$. We consider estimating $\boldsymbol{\beta}^*$ as follows

$$(2.3) \quad (\hat{b}_1, \ldots, \hat{b}_K, \widehat{\boldsymbol{\beta}}^{\mathrm{CQR}}) = \arg\min_{b_1,\ldots,b_k,\boldsymbol{\beta}} \sum_{k=1}^{K} \left\{ \sum_{i=1}^{n} \rho_{\tau_k}(y_i - b_k - \mathbf{x}_i^T \boldsymbol{\beta}) \right\}.$$

We call it *composite quantile regression*, for the objective function in (2.3) is a mixture of the objective functions from different quantile regression models. Typically, we use the equally spaced qauntiles: $\tau_k = \frac{k}{K+1}$ for $k = 1, 2, \ldots, K$.

We now establish the asymptotic normality of $\widehat{\boldsymbol{\beta}}^{\mathrm{CQR}}$. The following two regularity conditions are assumed throughout the rest of our discussions:

(1) There is a $p \times p$ positive definite matrix $\mathbf{C}$ such that

$$\lim_{n \to \infty} \frac{1}{n} \mathbf{X}^T \mathbf{X} = \mathbf{C}$$

where $\mathbf{C}$ is a $p \times p$ positive definite matrix.
(2) $\varepsilon$ has cumulative distribution function $F(\cdot)$ and density function $f(\cdot)$. For each $p$-vector $u$,

$$\lim_{n \to \infty} \frac{1}{n} \sum_{i=1}^{n} \int_0^{u_0 + \mathbf{x}_i^T \mathbf{u}} \sqrt{n}[F(a + t/\sqrt{n}) - F(a)] \, dt$$
$$= \frac{1}{2} f(a)(u_0, \mathbf{u}^T) \begin{bmatrix} 1 & 0 \\ 0 & \mathbf{C} \end{bmatrix} (u_0, \mathbf{u}^T)^T.$$

Conditions (1)–(2) are basically the same conditions for establishing the asymptotic normality of a single quantile regression [Koenker (2005)]. Under these conditions, we have the following result for the CQR estimate.

THEOREM 2.1 (The limiting distribution). *Under the regularity conditions* (1) *and* (2), *the limiting distribution of* $\sqrt{n}(\widehat{\boldsymbol{\beta}}^{\mathrm{CQR}} - \boldsymbol{\beta}^*)$ *is* $N(0, \boldsymbol{\Sigma}_{\mathrm{CQR}})$ *where*

$$\boldsymbol{\Sigma}_{\mathrm{CQR}} = \mathbf{C}^{-1} \frac{\sum_{k,k'=1}^{K} \min(\tau_k, \tau_{k'})(1 - \max(\tau_k, \tau_{k'}))}{(\sum_{k=1}^{K} f(b_{\tau_k}^*))^2}.$$

**3. Asymptotic relative efficiency.** In this section, we investigate the asymptotic relative efficiency (ARE) of the CQR with respect to the least squares. The same results can be applied to compute the relative efficiency of the CQR-oracle with respect to the LS-oracle.

Note that when $\sigma^2 < \infty$, the asymptotic variance of the least squares is $\sigma^2 \mathbf{C}^{-1}$. Therefore, the ARE of the CQR with respect to the least squares is

$$\text{ARE}(\tau_1, \ldots, \tau_K, f) = \frac{\sigma^2 (\sum_{k=1}^{K} f(b_{\tau_k}^*))^2}{\sum_{k,k'=1}^{K} \min(\tau_k, \tau_{k'})(1 - \max(\tau_k, \tau_{k'}))}. \quad (3.1)$$

We define the CQR-oracle estimator as follows

$$(\hat{b}_1, \ldots, \hat{b}_K, \widehat{\boldsymbol{\beta}}^{\text{CQR}}(oracle)_{\mathcal{A}})$$
$$= \underset{b_1, \ldots, b_K, \boldsymbol{\beta}}{\arg\min} \sum_{k=1}^{K} \left\{ \sum_{i=1}^{n} \rho_{\tau_k} \left( y_i - b_k - \sum_{j=1}^{q} x_{ij} \beta_j \right) \right\}, \quad (3.2)$$

and $\widehat{\boldsymbol{\beta}}^{\text{CQR}}(oracle)_{\mathcal{A}^c} = 0$. By Theorem 2.1 we have

$$\sqrt{n}(\widehat{\boldsymbol{\beta}}^{\text{CQR}}(oracle)_{\mathcal{A}} - \boldsymbol{\beta}_{\mathcal{A}}^*) \to_d N(0, \boldsymbol{\Sigma}_{\text{CQRoracle}}), \quad (3.3)$$

where

$$\boldsymbol{\Sigma}_{\text{CQRoracle}} = \mathbf{C}_{\mathcal{A}\mathcal{A}}^{-1} \frac{\sum_{k,k'=1}^{K} \min(\tau_k, \tau_{k'})(1 - \max(\tau_k, \tau_{k'}))}{(\sum_{k=1}^{K} f(b_{\tau_k}^*))^2}.$$

For the LS-oracle we have

$$\sqrt{n}(\widehat{\boldsymbol{\beta}}^{\text{LS}}(oracle)_{\mathcal{A}} - \boldsymbol{\beta}_{\mathcal{A}}^*) \to_d N(0, \sigma^2 \mathbf{C}_{\mathcal{A}\mathcal{A}}^{-1}). \quad (3.4)$$

Therefore, the asymptotic relative efficiency (ARE) of the CQR-oracle with respect to the LS-oracle is also equal to $\text{ARE}(\tau_1, \ldots, \tau_K, f)$ given in (3.1).

Take $\tau_k = \frac{1}{K+1}$ and write $\text{ARE}(K, f) \equiv \text{ARE}(\tau_1, \ldots, \tau_K, \sigma^2, f)$. It turns out that as $K$ approaches infinity $\text{ARE}(K, f)$ converges to a limit, denoted by $\delta(f)$. The next theorem gives us the explicit expression of $\delta(f)$ and provides a universal lower bound to $\delta(f)$.

THEOREM 3.1 The universal lower bound.

$$\lim_{K \to \infty} \frac{\sum_{k,k'=1}^{K} \min(\tau_k, \tau_{k'})(1 - \max(\tau_k, \tau_{k'}))}{(\sum_{k=1}^{K} f(b_{\tau_k}^*))^2} = \frac{1}{12(E_\varepsilon[f(\varepsilon)])^2}$$

and

$$\delta(f) \equiv \lim_{K \to \infty} \text{ARE}(K, f) = 12\sigma^2 (E_\varepsilon[f(\varepsilon)])^2.$$

Denote by $\mathcal{F}$ the collection of all density functions that satisfy condition (2) and have a finite variance. We have

$$\inf_{f \in \mathcal{F}} \delta(f) > \frac{6}{e\pi} = 0.7026.$$

Although $\delta(f)$ explicitly depends on $\sigma^2$, it is actually scale-invariant. We should also point out that the lower bound 70.26% given above is conservative. For commonly used error distributions in practice, $\delta$ is often much larger than the lower bound. Having the lower bound is a very useful property. It prevents severe loss in efficiency when using the CQR estimator instead of the LS estimator. Even in the worst possible scenario, the potential loss in efficiency is less than 30%. Meanwhile, CQR can have big gain in efficiency compared to the LS estimator, if the error follows certain types of error distributions as shown in the following examples.

We now calculate $\delta$ for some commonly used distributions.

EXAMPLE 1 (Normal distribution). Suppose the error density is $f(\varepsilon) = \frac{1}{\sqrt{2\pi}} e^{-\varepsilon^2/2}$. Then the least squares is the most efficient. We calculate

$$E_\varepsilon[f(\varepsilon)] = \int_{-\infty}^{\infty} \frac{1}{2\pi} e^{-\varepsilon^2} d\varepsilon = \frac{1}{2\sqrt{\pi}}.$$

Thus $\delta = \frac{3}{\pi} = 0.955$. In other words, the CQR is almost as efficient as the least squares in this case.

EXAMPLE 2 (Double exponential distribution). The density function is $f(\varepsilon) = \frac{1}{2} e^{-|\varepsilon|}$. We compute

$$E_\varepsilon[f(\varepsilon)] = \int_{-\infty}^{\infty} \tfrac{1}{4} e^{-2|\varepsilon|} d\varepsilon = \tfrac{1}{4}.$$

The error variance is 2. Hence Theorem 3.1 says $\delta = 2 \cdot 12 \cdot \frac{1}{16} = 1.5$.

EXAMPLE 3 (Logistic distribution). The density function of the logistic distribution is $f(\varepsilon) = \frac{e^\varepsilon}{(e^\varepsilon+1)^2}$. We compute

$$E_\varepsilon[f(\varepsilon)] = \int_{-\infty}^{\infty} \frac{e^{2\varepsilon}}{(1+e^\varepsilon)^4} d\varepsilon = \int_0^{\infty} \frac{s}{(1+s)^4} ds = \frac{1}{6}.$$

By Theorem 3.1 we know

$$\lim_{K \to \infty} \boldsymbol{\Sigma}_{\mathrm{CQR}} = 3\mathbf{C}^{-1}.$$

On the other hand, the Fisher information of the logistic distribution is $\frac{1}{3}$. Therefore the CQR can asymptotically achieve the information bound if the error follows the logistic distribution. Moreover, the variance of logistic distribution is $\frac{\pi^2}{3}$. Thus the relative efficiency is $\frac{\pi^2}{3} \cdot 12 \cdot \frac{1}{36} = \frac{\pi^2}{9} = 1.097$.

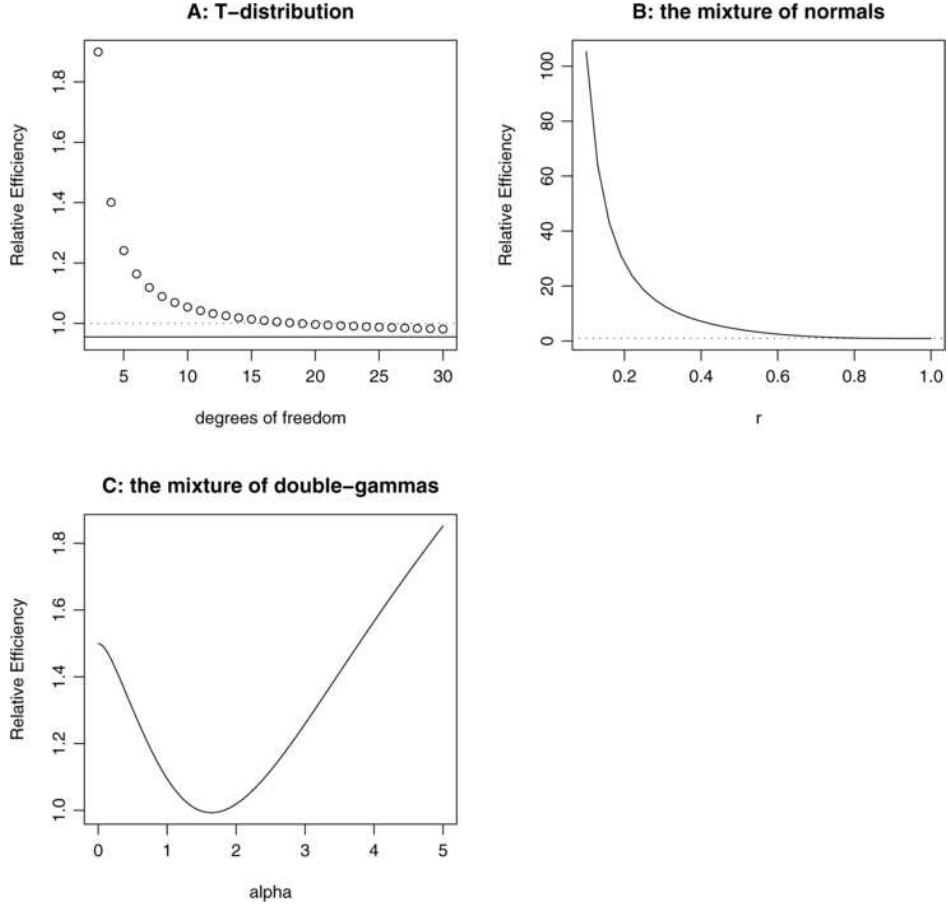

FIG. 1. (A): The relative efficiency ($\delta$) as a function of the degrees of freedom of the T-distribution. The dotted horizontal line indicates $\delta = 1$, while the solid line indicates $\delta = 0.955$. (B): The relative efficiency ($\delta$) as a function of $r$ in the mixture of normals distribution. The dotted horizontal line indicates $\delta = 1$. (C): The relative efficiency ($\delta$) as a function of $\alpha$ in the mixture of double gamma distribution.

EXAMPLE 4 (T-distribution). Let us also consider the T-distribution, which is often used to model errors following a heavy-tailed distribution.

COROLLARY 3.1. *For the T-distribution with degrees of freedom $v > 2$,*
$$\delta = \frac{12}{\pi} \frac{1}{v-2} \left( \frac{\Gamma((v+1)/2)}{\Gamma(v/2)} \right)^4 \left( \frac{\Gamma(v+1/2)}{\Gamma(v+1)} \right)^2.$$

For $v = 3$ $\delta$ is 1.9. From panel (A) in Figure 1, we see that the relative efficiency is greater than 1 for small degrees of freedoms. It is also interesting to see that for large degrees of freedom the relative efficiency is very close to 0.955. This is expected because the T-distribution converges to the normal as $v \to \infty$.

EXAMPLE 5 (A mixture of two normals). We have seen that for the normal distributed errors, $\delta$ is 0.955. It turns out that we can let $\delta$ be arbitrarily large by slightly perturbing the normal distribution, while keeping the error variance bounded.

COROLLARY 3.2. *Suppose the error follows a distribution of the mixture of normals*

$$\varepsilon \sim (1-r)N(0,1) + rN(0,r^6)$$

*for $0 < r < 1$. Then*

$$\delta = \frac{3}{\pi}\left((1-r)^2 + \frac{1}{r} + \frac{2\sqrt{2}r(1-r)}{\sqrt{1+r^6}}\right)^2 (1-r+r^7).$$

In Figure 1 panel (B) shows the curve of $\delta$ as a function of $r$. When $r$ is close to 0, the error variance approaches 1, but $\delta \approx \frac{3}{\pi}(1+\frac{1}{r})^2 \to \infty$.

EXAMPLE 6 (A mixture of two double Gamma distributions). We say $\varepsilon$ follows a double Gamma distribution with parameters $\alpha$ if

$$f(\varepsilon) = \frac{1}{2}\frac{1}{\Gamma(\alpha+1)}|\varepsilon|^\alpha e^{-|\varepsilon|}.$$

The double exponential distribution is a special double Gamma distribution using $\alpha = 0$.

COROLLARY 3.3. *Consider a mixture of double gamma distributions as follows:*

$$\varepsilon \sim e^{-\alpha}\frac{1}{2}e^{-|\varepsilon|} + (1-e^{-a})\frac{1}{\Gamma(\alpha+1)}\varepsilon^\alpha e^{-|\varepsilon|},$$

*where $\alpha \geq 0$. Then $\delta$ is equal to*

$$12(2e^{-\alpha} + (1-e^{-\alpha})(\alpha+1)(\alpha+2))$$
$$\times \left(\frac{e^{-2\alpha}}{4} + \frac{e^{-\alpha}(1-e^{-\alpha})}{2^{\alpha+1}} + \frac{(1-e^{-\alpha})^2\Gamma(2\alpha+1)}{4^{\alpha+1}\Gamma^2(\alpha+1)}\right)^2.$$

Using Stirling's formula [Feller (1968)] we have $\frac{\Gamma(2\alpha+1)}{4^{\alpha+1}\Gamma^2(\alpha+1)} \approx \frac{1}{4\sqrt{\pi\alpha}}$. Thus we can show that for $\alpha \to \infty$, $\delta \approx \frac{3}{4\pi}\alpha$. Displayed in panel (C) of Figure 1 is the curve of $\delta$ as a function of $\alpha$.

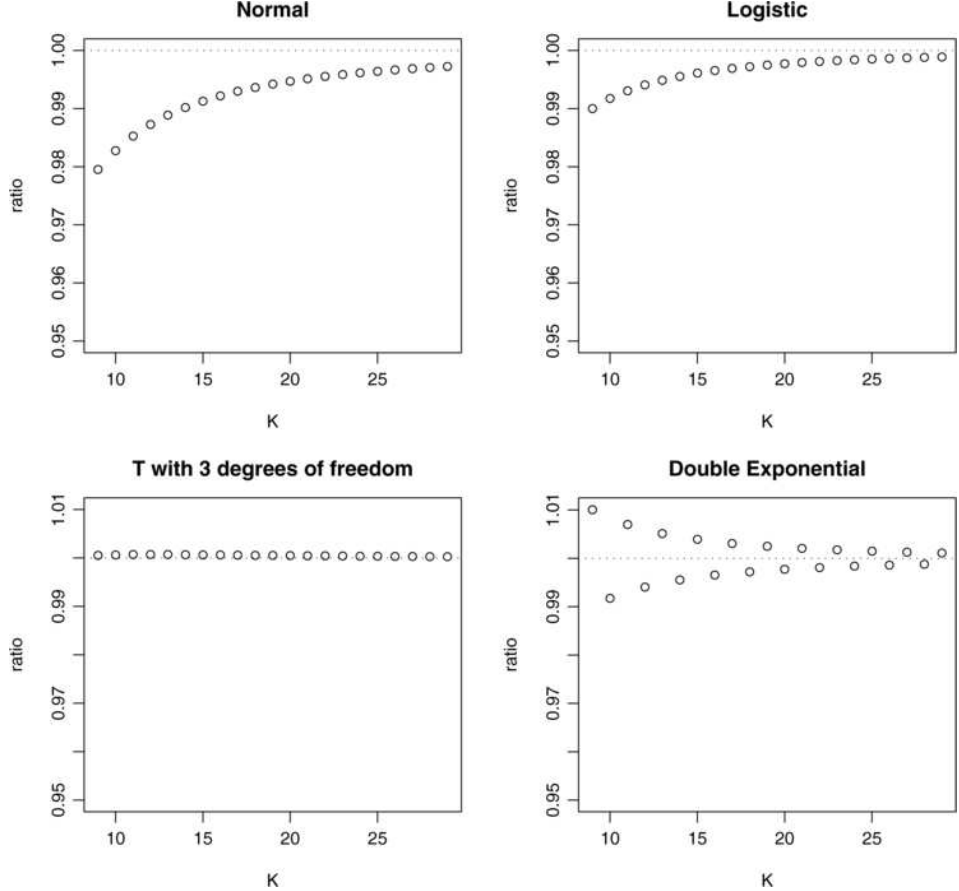

FIG. 2. *The approximation accuracy is measured by the ratio $\frac{\mathrm{RE}(K,f)}{12\sigma^2(E_\varepsilon[f(\varepsilon)])^2}$. We see that the ratio is almost 1 for each $K$ in $9, 10, \ldots, 29$. The dotted line in each panel indicates 1.*

In the above discussion we have considered the limit of the relative efficiency when $K \to \infty$. Empirically, we have found that for a reasonably large $K$, $\mathrm{RE}(K,f)$ is already very close to its limit. The ratio $\frac{\mathrm{RE}(K,f)}{12\sigma^2(E_\varepsilon[f(\varepsilon)])^2}$ measures the approximation accuracy. As can be seen from Figure 2, the ratio is very close to 1 for $K \geq 9$ in all the four different distributions considered there. In practice, it seems that $K = 19$ is a good choice, which amounts to using the $5\%, 10\%, 15\%, \ldots, 95\%$ quantiles.

**4. The CQR-oracular estimator.** The oracle model selection theory of Fan and Li (2001) contains two parts. The first part defines an optimal oracle estimator and the second part creates a practical procedure to achieve the optimal properties of the oracle. Following Fan and Li (2001), we say an estimation procedure $\xi$ is a *CQR-oracular estimator*, if $\widehat{\boldsymbol{\beta}}(\xi)$ (asymptotically) has the following two properties:

- Consistent selection: $\Pr(\{j : \widehat{\boldsymbol{\beta}}(\xi)_j \neq 0\} = \mathcal{A}) \to 1$.
- Efficient estimation: $\sqrt{n}(\widehat{\boldsymbol{\beta}}(\xi)_\mathcal{A} - \boldsymbol{\beta}^*_\mathcal{A}) \to_d N(0, \boldsymbol{\Sigma}_{\text{CQRoracle}})$.

Adaptive penalization methods have been successfully used to produce LS-oracular estimators. Fan and Li (2001) proposed the SCAD-penalized least squares and proved its oracle properties. Zou (2006) proposed the adaptive lasso and proved its oracle properties. In this section we present a CQR-oracular estimator and prove its oracle properties.

We adopt the adaptive lasso idea from Zou (2006). Suppose we first fit the CQR estimator using all the predictors. Theorem 2.1 says that $\widehat{\boldsymbol{\beta}}^{\text{CQR}}$ is root-$n$ consistent. Then we use $\widehat{\boldsymbol{\beta}}^{\text{CQR}}$ to construct the adaptively weighted lasso penalty and consider the penalized CQR estimator as follows:

$$
\begin{aligned}
(4.1) \quad & (\hat{b}_1, \ldots, \hat{b}_K, \widehat{\boldsymbol{\beta}}^{\text{ACQR}}) \\
& = \underset{b_1, \ldots, b_k, \boldsymbol{\beta}}{\arg\min} \sum_{k=1}^{K} \left\{ \sum_{i=1}^{n} \rho_{\tau_k}(y_i - b_k - \mathbf{x}_i^T \boldsymbol{\beta}) \right\} + \lambda \sum_{j=1}^{p} \frac{|\beta_j|}{|\hat{\beta}_j^{\text{CQR}}|^2}.
\end{aligned}
$$

We show that the adaptive lasso penalized CQR estimator (ACQR) enjoys the oracle properties of the CQR-oracle.

THEOREM 4.1 (Oracle properties). *Assume the two regularity conditions in Theorem 2.1. If $\frac{\lambda}{\sqrt{n}} \to 0$ and $\lambda \to \infty$, then $\widehat{\boldsymbol{\beta}}^{\text{ACQR}}$ must satisfy*

1. *Consistency in selection*: $\Pr(\{j : \widehat{\beta}_j^{\text{ACQR}} \neq 0\} = \mathcal{A}) \to 1$.
2. *Asymptotic normality*: $\sqrt{n}(\widehat{\boldsymbol{\beta}}_\mathcal{A}^{\text{ACQR}} - \boldsymbol{\beta}^*_\mathcal{A}) \to N(0, \boldsymbol{\Sigma}_{\text{CQRoralce}})$.

We have two remarks.

1. The results in Section 3 can be directly applied to compare the efficiency of the ACQR with any LS-oracular estimator. Let $\eta$ be a LS-oracular estimator. Then $\sqrt{n}(\widehat{\boldsymbol{\beta}}_\mathcal{A}^{\eta} - \boldsymbol{\beta}^*_\mathcal{A}) \to N(0, \mathbf{C}_{\mathcal{A}\mathcal{A}}^{-1}\sigma^2)$, if $\sigma^2 < \infty$. The relative efficiency of the ACQR with respect to $\eta$ is $\text{ARE}(\tau_1, \tau_2, \ldots, \tau_K, f)$ in (3.1). Therefore, the relative efficiency of the ACQR compared to $\eta$ is always larger than 0.70 and can greatly exceed 1 for some error distributions.
2. We can also use the SCAD penalty in the (4.1) and the resulting estimator should also possess the oracle properties of the CQR-oracle. We choose the adaptive lasso penalty only for the computational consideration. Note that similar to ordinary quantile regression, computing $\widehat{\boldsymbol{\beta}}^{\text{ACQR}}$ is equivalent to a linear programming problem. Thus we can efficiently compute the ACQR estimator using the standard linear program solver.

**5. A simulation study.** In this section we use simulation to compare the LS-oracle and the CQR-oracle and examine the performance of the ACQR with finite samples. Our simulated data consist of a training set and an independent validation set. Models were fitted on training data only, and the validation data were used to select the tuning parameters. We simulated 100 data consisting of 100 training observations and 100 validation observations from the model

$$y = \mathbf{x}^T \boldsymbol{\beta} + \varepsilon,$$

where $\boldsymbol{\beta} = (3, 1.5, 0, 0, 2, 0, 0, 0)$ and the predictors $(x_1, x_2, x_3, x_4, x_5, x_6, x_7, x_8)$ follow a multivariate normal distribution $N(0, \boldsymbol{\Sigma_x})$ with $(\boldsymbol{\Sigma_x})_{i,j} = 0.5^{|i-j|}$ for $1 \le i, j \le 8$. This regression model was considered in Tibshirani (1996) and Fan and Li (2001). Here we considered five different error distributions.

EXAMPLE 1. $\varepsilon \sim N(0, 3)$.

EXAMPLE 2. $\varepsilon = \sigma \varepsilon^*$ where $\varepsilon^*$ follows the mixture of normal distribution as in Corollary 3.2 with $r = 0.5$. We let $\sigma = \sqrt{6}$.

EXAMPLE 3. $\varepsilon = \sigma \varepsilon^*$ where $\varepsilon^*$ follows the mixture of normal double gamma as in Corollary 3.3 with $\alpha = 14$. We let $\sigma = \frac{1}{9}$.

EXAMPLE 4. The error distribution is T-distribution with 3 degrees of freedom.

EXAMPLE 5. The error distribution is Cauchy.

We used the quantiles $\tau_k = \frac{k}{20}$ for $k = 1, 2, \ldots, 19$ in the CQR-oracle and the ACQR. The model error is computed by

$$\mathrm{ME} = E[(\widehat{\boldsymbol{\beta}} - \boldsymbol{\beta})^T \boldsymbol{\Sigma_x} (\widehat{\boldsymbol{\beta}} - \boldsymbol{\beta}^T)].$$

We use the notation $(NC, NIC)$ to denote the variable selection result, where $NC$ denotes the number of predictors in $\{x_1, x_2, x_5\}$ that have nonzero coefficient vectors, and $NIC$ denotes the number of predictors in $\{x_3, x_4, x_6, x_7, x_8\}$ that have nonzero coefficient vectors. Table 1 shows the average model errors and variable selection results over 100 replications. In the asymptotic sense, the LS-oracle is the best in Example 1, while the CQR-oracle works better in Examples 2–4. The numerical experiments agree with the theory. We also see that the model error of the ACQR is close to that of the CQR-oracle. Table 1 shows that the ACQR does an excellent job in variable selection.

Example 5 is different from Examples 1–4, because the error distribution has infinite variance in Example 5. The LS-oracle is not the optimal estimator in this case. The simulation confirmed the theory. The model error of the LS-oracle is more than 2500. The CQR-oracle and the ACQR still work very well in Example 5.

TABLE 1
*Simulation results*

|  |  | Example 1 | Example 2 | Example 3 | Example 4 | Example 5 |
|---|---|---|---|---|---|---|
|  | LS-oracle | 0.079 | 0.104 | 0.091 | 0.082 | 2788 |
| Model error | CQR-oracle | 0.085 | 0.033 | 0.043 | 0.060 | 0.134 |
|  | ACQR | 0.112 | 0.046 | 0.048 | 0.077 | 0.174 |
| Variable selection | ACQR | (3, 0.53) | (3, 0.21) | (3, 0.23) | (3, 0.39) | (3, 0.53) |

## 6. Proofs.

PROOF OF THEOREM 2.1.
Let $\sqrt{n}(\widehat{\boldsymbol{\beta}}^{\mathrm{CQR}} - \boldsymbol{\beta}^*) = \mathbf{u}_n$ and $\sqrt{n}(\hat{b}_k - b^*_{\tau_k}) = u_{n,k}$. Then $(u_{n,1}, \ldots, u_{n,K}, \mathbf{u}_n)$ is the minimizer of the following criterion:

$$L_n = \sum_{k=1}^{K} \sum_{i=1}^{n} \left( \rho_{\tau_k}\left(\varepsilon_i - b^*_{\tau_k} - \frac{u_k + \mathbf{x}_i^T \mathbf{u}}{\sqrt{n}}\right) - \rho_{\tau_k}(\varepsilon_i - b^*_{\tau_k}) \right).$$

By the identity [Knight (1998)]

$$|r - s| - |r| = -s(I(r > 0) - I(r < 0)) + 2\int_0^s [I(r \le t) - I(r \le 0)] \, dt,$$

we have

$$\rho_\tau(r - s) - \rho_\tau(r) = s(I(r < 0) - \tau) + \int_0^y [I(r \le t) - I(r \le 0)] \, dt.$$

Thus we write $L_n$ as follows:

$$L_n = \sum_{k=1}^{K} \sum_{i=1}^{n} \frac{u_k + \mathbf{x}_i^T \mathbf{u}}{\sqrt{n}} (I(\varepsilon_i < b^*_{\tau_k}) - \tau_k)$$
$$+ \sum_{k=1}^{K} \sum_{i=1}^{n} \int_0^{(u_k + \mathbf{x}_i^T \mathbf{u})/\sqrt{n}} [I(\varepsilon_i \le b^*_{\tau_k} + t) - I(\varepsilon_i \le b^*_{\tau_k})] \, dt$$
$$= \sum_{k=1}^{K} z_{n,k} u_k + \mathbf{z}_n^T \mathbf{u} + \sum_{k=1}^{K} B_n^{(k)},$$

where

$$z_{n,k} \equiv \frac{1}{\sqrt{n}} \sum_{i=1}^{n} (I(\varepsilon_i < b^*_{\tau_k}) - \tau_k),$$

$$\mathbf{z}_n \equiv \frac{1}{\sqrt{n}} \sum_{i=1}^{n} \mathbf{x}_i^T \left[ \sum_{k=1}^{K} (I(\varepsilon_i < b^*_{\tau_k}) - \tau_k) \right],$$

$$B_n^{(k)} \equiv \sum_{i=1}^n \int_0^{(u_k+\mathbf{x}_i^T\mathbf{u})/\sqrt{n}} [I(\varepsilon_i \leq b_{\tau_k}^* + t) - I(\varepsilon_i \leq b_{\tau_k}^*)]\, dt.$$

By the Cramér–Wald device and CLT, we know

$$(z_{n,1},\ldots,z_{n_k},\mathbf{z}_n^T)^T \to_d (z_1,\ldots,z_k,\mathbf{z}^T)^T \sim N(0,\boldsymbol{\Sigma}) \qquad \text{for some } \boldsymbol{\Sigma}$$

and

$$\sum_{k=1}^K z_{n,k} u_k + \mathbf{z}_n^T \mathbf{u} \to_d \sum_{k=1}^K z_k u_k + \mathbf{z}^T \mathbf{u}.$$

Moreover, we have

$$E[B_n^{(k)}] = \sum_{i=1}^n \int_0^{(u_k+\mathbf{x}_i^T\mathbf{u})/\sqrt{n}} [F(t + b_{\tau_k}^*) - F(b_{\tau_k}^*)]\, dt$$

$$= \frac{1}{n}\sum_{i=1}^n \int_0^{u_k+\mathbf{x}_i^T\mathbf{u}} \sqrt{n}\left[F\left(\frac{t}{\sqrt{n}} + b_{\tau_k}^*\right) - F(b_{\tau_k}^*)\right] dt$$

$$\to \frac{1}{2} f(b_{\tau_k}^*)(u_k, \mathbf{u}^T) \begin{bmatrix} 1 & 0 \\ 0 & \mathbf{C} \end{bmatrix} (u_k, \mathbf{u}^T)^T.$$

$$\mathrm{Var}[B_n^{(k)}]$$
$$= \sum_{i=1}^n E\left(\int_0^{(u_k+\mathbf{x}_i^T\mathbf{u})/\sqrt{n}} (I(\varepsilon_i \leq b_{\tau_k}^* + t) - I(\varepsilon_i \leq b_{\tau_k}^*)\right.$$
$$\left. - [F(b_{\tau_k}^* + t) - F(b_{\tau_k}^*)])\, dt\right)^2$$
$$\leq \sum_{i=1}^n E\left[\left|\int_0^{(u_k+\mathbf{x}_i^T\mathbf{u})/\sqrt{n}} (I(\varepsilon_i \leq b_{\tau_k}^* + t) - I(\varepsilon_i \leq b_{\tau_k}^*)\right.\right.$$
$$\left.\left. - [F(b_{\tau_k}^* + t) - F(b_{\tau_k}^*)])\, dt\right|\right]$$
$$\times 2\left|\frac{u_k+\mathbf{x}_i^T\mathbf{u}}{\sqrt{n}}\right|$$
$$\leq 4E[B_n^{(k)}]\frac{\max_{1\leq i \leq n}|u_k+\mathbf{x}_i^T\mathbf{u}|}{\sqrt{n}}$$
$$\to 0.$$

Hence

$$B_n^{(k)} \to_p \frac{1}{2} f(b_{\tau_k}^*)(u_k, \mathbf{u}^T) \begin{bmatrix} 1 & 0 \\ 0 & \mathbf{C} \end{bmatrix} (u_k, \mathbf{u}^T)^T.$$

Thus it follows that

$$L_n \to_d V(u_1, \ldots, u_k, \mathbf{u})$$
$$= \sum_{k=1}^{K} z_{n,k} u_k + \mathbf{z}^T \mathbf{u} + \sum_{k=1}^{K} \tfrac{1}{2} f(b^*_{\tau_k}) u_k^2 + \tfrac{1}{2} \left(\sum_{k=1}^{K} f(b^*_{\tau_k})\right) \mathbf{u}^T \mathbf{C} \mathbf{u}.$$

Since $L_n$ is a convex function, then following Knight (1998) and Koenker (2005) we have

$$\hat{\mathbf{u}}_n \to_d \left[\left(\sum_{k=1}^{K} f(b^*_{\tau_k})\right) \mathbf{C}\right]^{-1} \mathbf{z} \sim N\left(0, \left(\sum_{k=1}^{K} f(b^*_{\tau_k})\right)^{-2} \mathbf{C}^{-1} \mathbf{\Sigma}_{\mathbf{z}} \mathbf{C}^{-1}\right).$$

However,

$$\mathbf{\Sigma}_{\mathbf{z}} = \mathbf{C} \operatorname{Var}\left(\sum_{k=1}^{K} [I(\varepsilon < b^*_{\tau_k}) - \tau_k]\right) = \mathbf{C}\left[\sum_{k,k'=1}^{K} \min(\tau_k, \tau_{k'})(1 - \max(\tau_k, \tau_{k'}))\right].$$

Therefore,

$$\sqrt{n}(\widehat{\boldsymbol{\beta}}^{\mathrm{CQR}} - \boldsymbol{\beta}^*) \to_d N\left(0, \mathbf{C}^{-1} \frac{\sum_{k,k'=1}^{K} \min(\tau_k, \tau_{k'})(1 - \max(\tau_k, \tau_{k'}))}{(\sum_{k=1}^{K} f(b^*_{\tau_k}))^2}\right). \quad \Box$$

PROOF OF THEOREM 3.1. First, it is easy to check that if $\tau_k = \frac{k}{K+1}$, then

$$\frac{1}{K^2} \sum_{k,k'=1}^{K} \min(\tau_k, \tau_{k'})(1 - \max(\tau_k, \tau_{k'})) \to \frac{1}{12}.$$

On the other hand,

$$\frac{1}{K} \sum_{k=1}^{K} f(b^*_{\tau_k}) = \frac{1}{K} \sum_{k=1}^{K} f\left(F^{-1}\left(\frac{k}{K+1}\right)\right)$$
$$\to \int_0^1 f(F^{-1}(s))\, ds$$
$$= E_U[f(F^{-1}(U))],$$

where $U \sim \mathrm{Unif}(0,1)$. Note that $F^{-1}(U)$ follows the distribution of $\varepsilon$, thus

$$E_U[f(F^{-1}(U))] = E_\varepsilon[f(\varepsilon)].$$

To prove the lower bound, first we use Jensen's inequality

$$E_\varepsilon[f(\varepsilon)] > \exp(E_\varepsilon[\log(f(\varepsilon))]).$$

Let $g(\varepsilon) = \frac{1}{\sqrt{2\pi\sigma^2}} e^{-(\varepsilon-\mu)^2/(2\sigma^2)}$, where $\mu$ is the mean of the error distribution. Then by the entropy inequality we have

$$E_\varepsilon\left[\log\left(\frac{f(\varepsilon)}{g(\varepsilon)}\right)\right] \geq 0.$$

On the other hand, we compute

$$E_\varepsilon[\log(g(\varepsilon))] = E_\varepsilon\left[\log\left(\frac{1}{\sqrt{2\pi\sigma^2}}\right) - \frac{(\varepsilon-\mu)^2}{2\sigma^2}\right]$$
$$= \log\left(\frac{1}{\sqrt{2\pi\sigma^2}}\right) - \frac{1}{2}.$$

Thus we have

$$E_\varepsilon[f(\varepsilon)] > \exp\left(\log\left(\frac{1}{\sqrt{2\pi\sigma^2}}\right) - \frac{1}{2}\right) = \frac{1}{\sqrt{2e\pi\sigma^2}}.$$
$$\delta > 12\sigma^2\left(\frac{1}{\sqrt{2e\pi\sigma^2}}\right)^2 = \frac{6}{\pi e} = 0.7026. \qquad \square$$

PROOFS OF COROLLARIES 3.1–3.3. The proofs are direct applications of Theorem 3.1. The details are given in a technical report of the paper [Zou and Yuan (2007)], thus omitted for the sake of space. $\square$

PROOF OF THEOREM 4.1. We write $\lambda = \lambda_n$. Let $\sqrt{n}(\widehat{\boldsymbol{\beta}}^{\text{ACQR}} - \boldsymbol{\beta}^*) = \mathbf{u}_n$ and $\sqrt{n}(\hat{b}_k - b^*_{\tau_k}) = u_{n,k}$. Then $(u_{n,1}, \ldots, u_{n,K}, \mathbf{u}_n)$ is the minimizer of the following criterion:

$$L_n = \sum_{k=1}^{K}\sum_{i=1}^{n}\left(\rho_{\tau_k}\left(\varepsilon_i - b^*_{\tau_k} - \frac{u_k + \mathbf{x}_i^T\mathbf{u}}{\sqrt{n}}\right) - \rho_{\tau_k}(\varepsilon_i - b^*_{\tau_k})\right)$$
$$+ \sum_{j=1}^{p}\frac{\lambda_n}{\sqrt{n}|\hat{\beta}_j^{\text{CQR}}|^2}\sqrt{n}\left[\left|\beta_j^* + \frac{u_j}{\sqrt{n}}\right| - |\beta_j^*|\right].$$

Following the arguments in the proof of Theorem 2.1, we write $L_n$ as follows:

$$L_n = \left(\sum_{k=1}^{K} z_{n,k} u_k + \mathbf{z}_n^T \mathbf{u}\right) + \left(\sum_{k=1}^{K} B_n^{(k)}\right)$$
$$+ \left(\sum_{j=1}^{p}\frac{\lambda_n}{\sqrt{n}|\hat{\beta}_j^{\text{CQR}}|^2}\sqrt{n}\left[\left|\beta_j^* + \frac{u_j}{\sqrt{n}}\right| - |\beta_j^*|\right]\right).$$

If $\beta_j^* \neq 0$, then $|\hat{\beta}_j^{\text{CQR}}|^2 \to_p |\beta_j^*|^2$, and $\sqrt{n}(|\beta_j^* + \frac{u_j}{\sqrt{n}}| - |\beta_j^*|) \to u_j \text{sgn}(\beta_j^*)$. By Slutsky's theorem, $\frac{\lambda_n}{\sqrt{n}|\hat{\beta}_j^{\text{CQR}}|^2}\sqrt{n}(|\beta_j^* + \frac{u_j}{\sqrt{n}}| - |\beta_j^*|) \to_p 0$. If $\beta_j^* = 0$, then

$\sqrt{n}(|\beta_j^* + \frac{u_j}{\sqrt{n}}| - |\beta_j^*|) = |u_j|$, and $\frac{\lambda_n}{\sqrt{n}|\hat{\beta}_j^{\mathrm{CQR}}|^2} \to_p \infty$. Therefore, we have

$$\frac{\lambda_n}{\sqrt{n}|\hat{\beta}_j^{\mathrm{CQR}}|^2}\sqrt{n}\bigg(\bigg|\beta_j^* + \frac{u_j}{\sqrt{n}}\bigg| - |\beta_j^*|\bigg)$$

$$\to_p W(\beta_j, u_j) = \begin{cases} 0, & \text{if } \beta_j^* \neq 0, \\ 0, & \text{if } \beta_j^* = 0 \text{ and } u_j = 0, \\ \infty, & \text{if } \beta_j^* = 0 \text{ and } u_j \neq 0. \end{cases}$$

Thus it follows that

$$L_n \to_d \sum_{k=1}^K z_k u_k + \mathbf{z}^T \mathbf{u} + \sum_{k=1}^K \tfrac{1}{2} f(b_{\tau_k}^*) u_k^2$$

$$+ \tfrac{1}{2}\bigg(\sum_{k=1}^K f(b_{\tau_k}^*)\bigg)\mathbf{u}^T \mathbf{C} \mathbf{u} + \sum_{j=1}^p W(\beta_j, u_j).$$

Let us write $\mathbf{u} = (\mathbf{u}_1^T, \mathbf{u}_2^T)^T$ where $\mathbf{u}_1$ contains the first $q$ elements of $\mathbf{u}$. Using the same arguments in Knight (1998) and Koenker (2005), we have

$$\hat{\mathbf{u}}_{2,n} \to_d 0$$

and

$$\hat{\mathbf{u}}_{1,n} \to_d \bigg[\bigg(\sum_{k=1}^K f(b_{\tau_k}^*)\bigg)\mathbf{C}\bigg]^{-1} \mathbf{z} \sim N\bigg(0, \bigg(\sum_{k=1}^K f(b_{\tau_k}^*)\bigg)^{-2} \mathbf{C}_{\mathcal{A}\mathcal{A}}^{-1} \mathbf{\Sigma}_{\mathbf{z}_1} \mathbf{C}_{\mathcal{A}\mathcal{A}}^{-1}\bigg),$$

$$\mathbf{\Sigma}_{\mathbf{z}_1} = \mathbf{C}_{\mathcal{A}\mathcal{A}}\bigg[\sum_{k,k'=1}^K \min(\tau_k, \tau_{k'})(1 - \max(\tau_k, \tau_{k'}))\bigg].$$

Therefore, the asymptotic normality is proven.

We now prove the consistent selection result. Let $\widehat{\mathcal{A}}_n = \{j : \widehat{\boldsymbol{\beta}}_j^{\mathrm{ACQR}} \neq 0\}$. $\forall j \in \mathcal{A}$, the asymptotic normality indicates $\Pr(j \in \widehat{\mathcal{A}}_n) \to 1$. Then it suffices to show that $\forall j \notin \mathcal{A}$, $\Pr(j \in \widehat{\mathcal{A}}_n) \to 0$. We know $|\frac{\rho_\tau(r_1) - \rho_\tau(r_2)}{r_1 - r_2}| \leq \max(\tau, 1-\tau) < 1$. If $j \in \widehat{\mathcal{A}}_n$, then we must have $\frac{\lambda_n}{|\hat{\beta}_j^{\mathrm{CQR}}|^2} < \sum_{i=1}^n |x_{ij}|$. Thus we have $\Pr(j \in \widehat{\mathcal{A}}_n) \leq \Pr(\frac{\lambda_n}{|\hat{\beta}_j^{\mathrm{CQR}}|^2} < \sum_{i=1}^n |x_{ij}|)$. But $\frac{1}{n}\sum_{i=1}^n |x_{ij}| \leq \sqrt{\frac{1}{n}\sum_{i=1}^n |x_{ij}^2|} \to \mathbf{C}_{jj}$ and $\frac{\lambda_n}{n|\hat{\beta}^{\mathrm{CQR}}|^2} = \frac{\lambda_n}{|\sqrt{n}\hat{\beta}_j^{\mathrm{CQR}}|^2} \to \infty$, thus $P(j \in \widehat{\mathcal{A}}_n) \to 0$. $\square$

**7. Concluding remarks.** Fan and Li (2001) introduced the concept of oracle model estimator and proposed the SCAD method to achieve the oracle properties. Fan and Li (2001) showed that ideally if one knows the likelihood model, then the oracle is the maximum likelihood estimate knowing the true underlying sparse model, and the SCAD estimator is obtained by the penalized likelihood model using the SCAD penalty. However, in the linear regression problems, the error distribution (hence the likelihood model) is typically unknown. Hence we can only consider a practical oracle procedure. Fan and Li (2001) showed that penalized least squares using the SCAD penalty is a practical oracular estimator. Unfortunately, the LS-oracle and the SCAD procedure break down when the error distribution has an infinite variance. The other oracular method, the adaptive lasso [Zou (2006)], has the same trouble, because it also mimics the LS-oracle.

In this work we have proposed the composite quantile regression and proven its nice theoretical properties. We have shown that the oracle model selection theory of Fan and Li (2001) still works beautifully even for the cases where the error variance is infinite, as long as we replace the LS-oracle with the CQR-oracle. Compared with the LS-oracle, the CQR-oracle has two remarkable nice properties:

(1) Its relative efficiency is always larger than 70%.
(2) In the Gaussian model the relative efficiency of the ACQR is 95.5%. With nonnormal errors, its relative efficiency could be arbitrarily large.

Following the lines of Fan and Li (2001) and Zou (2006), we have developed the adaptively penalized CQR method (ACQR) and proven the oracle properties of the ACQR. There are a family of penalty functions, including the SCAD penalty, that can be used to create CQR-oracular estimators. We have used the adaptive lasso penalty only for its computational convenience.

As pointed out by the associate editor, we might also consider a general composite quantile regression problem by minimizing

$$(7.1) \qquad \int_0^1 \sum_{i=1}^n \rho_t(y_i - b_t - \mathbf{x}_i^T \boldsymbol{\beta}) w(t) \, dt,$$

where the weight function $w(t)$ is a density function over $(0, 1)$. The proposed CQR criterion uses a discrete uniform distribution on $\{\frac{1}{K+1}, \ldots, \frac{K}{K+1}\}$. Although the weight function $w(t)$ could be a continuous density function in (7.1), it seems that we need to discretize the integral in order to numerically compute the estimator. Hence, technically speaking, a discrete distribution density is used to construct the weights. In this paper we have shown that the discrete uniform distribution leads to an interesting estimator which enjoys various nice properties. It would be interesting to see whether other distributions could result in similar estimators. This is an open problem for future research.

**Acknowledgments.** The authors gratefully acknowledge the helpful comments of the Associate Editor and referees.

SCHOOL OF STATISTICS
UNIVERSITY OF MINNESOTA
MINNEAPOLIS, MINNESOTA 55455
USA
E-MAIL: hzou@stat.umn.edu

SCHOOL OF INDUSTRIAL AND
  SYSTEMS ENGINEERING
GEORGIA INSTITUTE OF TECHNOLOGY
ATLANTA, GEORGIA 30332-0205
USA
E-MAIL: myuan@isye.gatech.edu